\documentclass{article}
\usepackage{spconf,amsmath,graphicx}
\usepackage{amsthm}
\usepackage{amssymb}
\usepackage{mathtools}
\usepackage{tikz-cd}
\usepackage{pgfplots} 
\usepackage{pgf}
\usepackage{tikz}
\usepackage{enumitem}
\usepackage{algorithm}
\usepackage{enumitem}
\usepackage{algpseudocode}

\newtheorem{assumption}{Assumption}
\newtheorem{theorem}{Theorem}
\newtheorem{lemma}{Lemma}

\newcommand{\subscript}[2]{$#1 _ #2$}
\title{Eco-PANDA: A computationally economic, geometrically converging dual optimization method on time-varying undirected graphs}
\name{Marie Maros and Joakim Jald\'en\thanks{This project has received funding from the European Research Council (ERC) under the European Union's Horizon 2020 research and innovation programme (Grant Agreement no. 742648)}}
\address{Department of Information Science and Engineering, EECS \\
Royal Institute of Technology (KTH) \\
SE-100 44 Stockholm, Sweden}
\begin{document}
\ninept
\maketitle
\begin{abstract}
In this paper we consider distributed convex optimization over time-varying undirected graphs. We propose a linearized version of primarily averaged network dual ascent (PANDA) that keeps the advantages of PANDA while requiring less computational costs.
The proposed method, economic primarily averaged network dual ascent (Eco-PANDA),  provably converges at R-linear rate to the optimal point given that the agents' objective functions are strongly convex and have Lipschitz continuous gradients. Therefore, the method is competitive, in terms of type of rate, with both DIGing and PANDA. The proposed method halves the communication costs of methods like DIGing while still converging R-linearly and having the same per iterate complexity.
\end{abstract}
\section{Introduction}
\label{sec:intro}
In this paper we solve convex optimization problems of the form
\begin{equation}
\label{eq:original}
\underset{\mathbf{\bar{x}}\in\mathbb{R}^p}{\min} \quad \sum_{i=1}^n f_i(\mathbf{\bar{x}})
\end{equation}
in a decentralized manner within a network of $n$ agents. Agent $i$ possesses exclusive knowledge of its private objective function $f_i.$ Problems of the form of \eqref{eq:original} naturally arise in applications such as distributed estimation and control \cite{estimation1,estimation2,estimation3,estimation4}, decentralized source localization \cite{dsl}, power grids \cite{pg1,pg2} and distributed learning \cite{dl1,dl2,dl3}. Methods to solve problems of the form \eqref{eq:original} in a decentralized manner over static networks have received a lot of attention in the last few years. If the individual objective functions, $f_i,$ are strongly convex and have Lipschitz continuous gradients, linearly converging methods such as  the decentralized exact first-order algorithm (EXTRA) \cite{extra} have previously been proposed. 
In some applications, the agents may be communicating via wireless links. Due to the random nature of the wireless channel the communication links can be unreliable. This justifies the need for procedures according to which the agents cooperate to solve the optimization problem in \eqref{eq:original} over time varying networks. DIGing \cite{DIGing} is the first known method to have established linear convergence ver time-varying networks. However, DIGing requires the exchange of both the primal variables and the gradients at each iteration. 

Dual methods have been shown to reach optimal convergence rates in distributed optimization \cite{dual_optimal}. However, with the exception of PANDA \cite{cdc_panda,panda}, no dual method has been proposed to date that can provably converge linearly on time-varying networks. In \cite{panda} we establish PANDA's R-linear converge and demonstrate that it requires the exchange of only the primal variables as opposed to DIGing. However, at each iteration the nodes must solve a convex optimization problem. This is not problematic whenever the objective function is dual friendly \cite{dual_friendly}, but may be computationally expensive otherwise. To alleviate this we propose herein Eco-PANDA, a linearized version of PANDA, capable of converging linearly on time varying graphs while still requiring the exchange of only the primal variables at each iteration. 

The remainder of this paper is organized as follows. In Section \ref{sec:format} we formalize the problem further, introduce Eco-PANDA and provide the formal statement regarding Eco-PANDA's convergence. In Section \ref{sec:proof} we provide proof the statement made in Section \ref{sec:format}. In Section \ref{sec:typestyle} we empirically illustrate Eco-PANDA's convergence. Finally we end the paper with some concluding remarks.

\section{Problem Formulation and Algorithm}
\label{sec:format}

Let $\mathcal{G}(k) = (\mathcal{V},\mathcal{E}(k))$ denote the undirected graph representing the network that connects the $n$ agents at time $k.$ The agents are represented by the vertices $\mathcal{V}$ and the connections between the agents by the edges $\mathcal{E}(k).$ Agents $i$ and $j$ can communicate if the edge $(i,j) \in \mathcal{E}(k).$ Further, we assume the graph is undirected, and therefore, if $(i,j)\in \mathcal{E}(k)$ then $(j,i) \in \mathcal{E}(k).$ Further, the set of agents with which agent $i$ can communicate will be denoted $\mathcal{N}_i(k).$ In order to solve the optimization problem \eqref{eq:original} in a distributed manner over the graph $\mathcal{G}(k)$ we introduce for each of the agents a copy $\mathbf{x}_i$ of the optimization variable $\mathbf{\bar{x}}$ and define the optimization problem
\begin{subequations}
\label{eq:distributed}
\begin{align}
& \underset{\mathbf{x}\in \mathbb{R}^{np}}{\min} \quad f(\mathbf{x}) = \sum_{i=1}^n  f_i(\mathbf{x}_i) \\
& \text{s.t.} \quad \left(\mathbf{I}_{np}-\mathbf{W}(k) \otimes \mathbf{I}_p\right)\mathbf{x} = 0\,,
\end{align}
\end{subequations} 
where $\mathbf{x} \triangleq (\mathbf{x}_1^T,\hdots,\mathbf{x}_n^T)^T$ and $\mathbf{W}(k)$ fulfils the requirements provided in Assumption \ref{assumption:network}.
\begin{assumption}[Mixing Matrix Sequence $\{\mathbf{W}(k)\}$\label{assumption:network}]
For any $k=0,\hdots,1,$ the mixing matrix $\mathbf{W}(k) \in \mathbb{R}^{n \times n}$ satisfies the following relations:
\begin{enumerate}[label = (\subscript{P}{{\arabic*}})]
\item \label{consensus1} Decentralized property: if $i \neq j$ and $(i,j) \not\in \mathcal{E}(k)$ $w_{ij}(k) = 0,$ i.e., $\mathbf{W}(k)$ is defined on the edges of the graph $\mathcal{G}(k).$
\item \label{consensus2} Doubly stochastic: $\mathbf{W}(k) \mathbf{1}_n = \mathbf{1}_n,$ $\mathbf{1}_ n^T \mathbf{W}(k) = \mathbf{1}_n^T.$
\item \label{cycleB} Joint spectrum property: Let $\sigma_{\text{max}}$ denote the largest singular value of a matrix and let
\begin{equation}
\mathbf{W}_b(k) \triangleq \mathbf{W}(k)\mathbf{W}(k-1)\hdots \mathbf{W}(k-b+1),
\end{equation}
for $k \geq 0$ and $b \geq k-1,$ with $\mathbf{W}_b(k) = \mathbf{I}_n$ for $ k<0$ and $\mathbf{W}_0(k) = \mathbf{I}_n.$ Then, there exists a positive $B$ such that
\begin{equation}
\underset{k \geq B-1}{\sup} \quad \delta(k)=\delta <1,
\end{equation} 
where
\begin{equation}
\delta(k) = \sigma_{\text{max}}\left\{\mathbf{W}_B(k)-\frac{1}{n}\mathbf{1}_n^T\mathbf{1}_n\right\},\,\forall k=0,1,\hdots.
\end{equation}
\end{enumerate}
\end{assumption}
Properties \ref{consensus1} and \ref{consensus2} in Assumption \ref{assumption:network} are common in the consensus literature while \ref{cycleB} is a  requirement on the time-varying nature of the graph $\mathcal{G}(k).$ Examples and discussion on condition \ref{cycleB} can be found in \cite{DIGing}. 

Eco-PANDA follows the same intuition as PANDA. For an intuition on PANDA please refer to \cite{panda}. PANDA solves the optimization problem

\begin{equation}
\label{eq:fencheldual}
\underset{\mathbf{y}:((\mathbf{I}_n-\frac{1}{n}\mathbf{1}_n\mathbf{1}_n^T)\otimes \mathbf{I}_p)\mathbf{y}=\mathbf{0}}{\min} \quad f^{*}(\mathbf{y}),
\end{equation}
where $\mathbf{1}_n \in \mathbb{R}^n$ denotes the vector of all ones and $f^*(\mathbf{y}) \triangleq \underset{\mathbf{x} \in \mathbb{R}^{np}}{\min} \quad \mathbf{y}^T\mathbf{x}-f(\mathbf{x})$ denotes the Fenchel conjugate of $f$ by applying approximated projected gradient descent on \eqref{eq:fencheldual}, where the projection is approximated via a gradient tracking scheme much like the schemes in \cite{DIGing} and \cite{dynaconsensus}. While PANDA achieves R-linear convergence rates by communicating half as many variables as DIGing per iteration, the communication cost is decreased at the expense of increased computational costs. In particular, PANDA requires obtaining the minimizer
\begin{equation}
\label{eq:grad_fenchel}
\mathbf{x}_{\text{p}}(k+1):= \nabla f^*(\mathbf{y}(k))=\underset{\mathbf{x} \in \mathbb{R}^{np}}{\min} \quad f(\mathbf{x}) - \mathbf{y}(k)^T\mathbf{x},
\end{equation}
for a given value of $\mathbf{y}(k)$ at each iteration, while DIGing requires only a simple update in the direction of the approximated average gradient. Eco-PANDA addresses this by using a quadratic upper bound on the objective functions that leads to simple iterates.
Before introducing Eco-PANDA we will introduce the following assumption on the objective functions.
\begin{assumption}[Strong convexity and Lipschitz continuity of gradients\label{assumption:function}]
The objective function $f$ is $\mu-$strongly convex and its gradient is $L-$Lipschitz continuous.
\end{assumption}
In order to make PANDA competitive to DIGing in terms of iterate complexity  we propose to upper bound the objective function in \eqref{eq:grad_fenchel} by
\begin{equation}
\label{eq:upper bound}
f(\mathbf{x}(k)) + (\nabla f(\mathbf{x}(k))-\mathbf{y}(k))^T(\mathbf{x}-\mathbf{x}(k)) + \frac{\eta}{2}\|\mathbf{x}-\mathbf{x}(k)\|^2.
\end{equation}
Since the function $f$ is convex and is $L-$Lipschitz continuous, as long as $\eta > L$ \eqref{eq:upper bound} is an upper bound to the objective function in \eqref{eq:grad_fenchel}. Further, the upper bound becomes tight as $\mathbf{y}(k) \to \mathbf{y}^{\star}$ and $\mathbf{x}(k) \to \mathbf{x}^{\star}.$ With this we are ready to introduce Eco-PANDA which is formally written in Algorithm \ref{alg:lpanda}.

\begin{algorithm}
\caption{Eco-PANDA \label{alg:lpanda}}
\begin{algorithmic}[1]
\State Choose step size $c > 0,$ $\eta > L$ and pick $\mathbf{z}(0)= \mathbf{x}(0)$ and $\mathbf{y}(0)$ such that $(\boldsymbol{\Pi}_{\mathbf{1}_n} \otimes \mathbf{I}_p)\mathbf{y}(0) = \mathbf{0}.$
\For{$k=0,1,\hdots$} each agent $i=1,\hdots,n$:
\State computes $$
\mathbf{x}_i(k+1) :=  \mathbf{x}_i(k) -\frac{1}{\eta}(  \nabla f_i(\mathbf{x}_i) - \mathbf{y}_i(k))$$
\State exchanges $\mathbf{z}_i(k)$
 with the agents in $\mathcal{N}_i(k).$ 
\State computes 
$$\mathbf{z}_i(k+1) := \sum_{j \in \mathcal{N}_i(k) \cup \{i\}}w_{ij}(k)\mathbf{z}_j(k) + \mathbf{x}_i(k+1) - \mathbf{x}_i(k)$$
\State and computes $$\mathbf{y}_i(k+1) := \mathbf{y}_i(k) - c(\mathbf{x}_i(k+1) -\mathbf{z}_i(k+1))$$
\EndFor
\end{algorithmic}
\end{algorithm}

We are now ready to provide the paper's main statement.
\begin{theorem}[Eco-PANDA converges R-linearly\label{theorem:main}]
Let Assumptions \ref{assumption:network} and \ref{assumption:function} hold. Also, let $\kappa \triangleq \frac{L}{\mu}$ denote the condition number of $f$ an $q \triangleq 1-\frac{\mu}{\eta}.$
Define the quantities
\begin{align}
C \triangleq \frac{q(1-\sqrt{q})}{(\max\{\frac{1}{\eta},\frac{q}{\mu}\} + B(1-\delta))(3+q)}, \quad
\alpha \triangleq \frac{C}{\mu} < 1\\
\lambda_{\bar{c}}^B \triangleq \frac{\alpha\delta + \sqrt{4\kappa^{3/2}(4\kappa^{3/2}-\alpha\delta^2)}}{\alpha + \kappa^{3/2}},\,
\bar{c} \triangleq \frac{\alpha \mu (\lambda_{\bar{c}}^B - \delta)^2}{2\sqrt{\kappa}}
\end{align} Then for any step size
\begin{equation}
c \in \left(0, \frac{\alpha\mu(1-\delta)^2}{2}\right)
\end{equation}
the sequence $\{\mathbf{y}(k)\}_{k \geq 0}$ converges to $\mathbf{y}^{\star},$ the unique solution of \eqref{eq:fencheldual}, and the sequences $\{\mathbf{x}_i(k)\}_{k \geq 0}$ converge to $\mathbf{\bar{x}}^{\star}$ the unique solution of \eqref{eq:original}, at a global R-linear rate of at least $\mathcal{O}(\lambda^k),$ where $\lambda < 1$ is given by
\begin{equation}
\lambda = \begin{cases}
\sqrt[2B]{1-\frac{c}{2L}} & c \in (0,\bar{c}] \\
\sqrt[B]{\delta + \sqrt{\frac{2c\sqrt{\kappa}}{\alpha\mu}}} & c \in (\bar{c},\frac{\alpha \mu(1-\delta)^2}{2})
\end{cases}
\end{equation} 
\end{theorem}
When comparing the statement in Theorem \ref{theorem:main} to the main statement in \cite{panda} one can see that the upper bound on the convergence rate is affected by the parameter $\alpha.$ In particular, for PANDA $\alpha = 1,$ while for Eco-PANDA $\alpha < 1.$

\section{Proof of convergence}\label{sec:proof}
\begin{figure}
\begin{center}
\begin{tikzcd}[row sep=-0.35em, column sep = 0.4em]
&&&&& \mathbf{z}^{\perp} \ar{dr} && \\
\!\!\mathbf{r} \ar{r}
& \mathbf{x}^{\perp} \ar{r}
& \Delta_{\mathbf{xz}}^{\perp} \ar{r}
& \Delta\mathbf{y} \ar{r}
& \mathbf{s} \ar{ur}
                \ar{dr}
& 
& \boldsymbol{\epsilon} \ar{r}
& \mathbf{r} \\
&&&&& \Delta_{\mathbf{x}} \ar{ur} && 
\end{tikzcd},
\end{center}
\caption{\label{fig:circle}Cycle of arrows}
\end{figure}
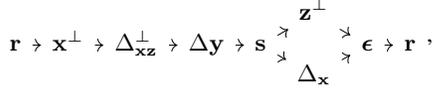

This section is devoted to providing the guidelines to proving Theorem \ref{theorem:main}. We will start the section by introducing some notation and the small gain theorem, which is the main tool used to establish the convergence of Eco-PANDA. Let $\mathbf{s}^i \triangleq \{\mathbf{s}^i(0),\mathbf{s}^i(1),\hdots\}$ denote an infinite sequence of vectors $\mathbf{s}^i(k) \in \mathbb{R}^{np},$ for $i=1,\hdots,m.$ Further, let
\begin{equation}
\|\mathbf{s}^i\|^{\lambda,K} \triangleq \underset{k=0,\hdots,K}{\sup} \quad \frac{1}{\lambda^k}\|\mathbf{s}^i(k)\|, \,\,
\|\mathbf{s}^i\|^{\lambda} \triangleq \underset{k \geq 0}{\sup} \quad \frac{1}{\lambda^k}\|\mathbf{s}^i(k)\|.
\end{equation}
\begin{theorem}[Small gain theorem \cite{DIGing}] Suppose $\mathbf{s}^1,\hdots,\mathbf{s}^m$ are vector sequences such that for all positive $k$ and for each $i=1,\hdots,m,$ we have an arrow $\mathbf{s}^i \to \mathbf{s}^{(i \mod m)+1},$ i.e.,
\begin{equation}
\label{eq:mod}
\|\mathbf{s}^{(i \mod m) + 1}\|^{\lambda,K} \leq \gamma_i\|\mathbf{s}^i\|^{\lambda,K} + \omega_i,
\end{equation} 
where the constants $\gamma_1,\hdots,\gamma_m$ and $\omega_1,\hdots,\omega_m$ are independent of $K.$ Further, suppose that the constants $\gamma_1,\hdots,\gamma_m$ are non-negative and satisfy
$\gamma_1\hdots\gamma_m < 1.$
Then we have that
\begin{align}
\|\mathbf{s}^1\|^{\lambda}& \leq \frac{1}{1-\gamma_1\gamma_2\hdots\gamma_m}(\omega_1\gamma_1\gamma_3\hdots\gamma_m + \omega_2\gamma_3\gamma_4\hdots\gamma_m \nonumber\\
& + \hdots + \omega_{m-1}\gamma_m + \omega_m).
\end{align}
Further, if $\|\mathbf{s}^1\|^{\lambda} < C$ where $C < \infty,$ $\|\mathbf{s}^1(k)\|$ converges to zero exponentially fast and at rate $\lambda.$ Proof of this statement can be found in \cite{DIGing}. Further, due to the cyclic nature of the small gain theorem all sequences will converge to zero exponentially fast at rate $\lambda.$ 
\end{theorem}
In order to establish the convergence of Eco-PANDA we will first establish the cycle of arrows in Fig. \ref{fig:circle}
where, given that $\boldsymbol{\Pi}_{\mathbf{1}_n}^{\perp} \triangleq (\mathbf{I}_n-\frac{1}{n}\mathbf{1}_n\mathbf{1}_n^T)\otimes \mathbf{I}_p,$ the sequences are defined as
 $\mathbf{r}(k) \triangleq \mathbf{y}(k) - \mathbf{y}^{\star},$ 
$\mathbf{x}^{\perp} \triangleq \boldsymbol{\Pi}_{\mathbf{1}_{n}}^{\perp}\mathbf{x}(k+1),$
$\Delta_{\mathbf{xz}}^{\perp}(k) \triangleq \boldsymbol{\Pi}_ {\mathbf{1}_n}^{\perp}(\mathbf{x}(k)-\mathbf{z}(k)),$
$\Delta\mathbf{y}(k) \triangleq \mathbf{y}(k) - \mathbf{y}(k-1),$
$\mathbf{s}(k) \triangleq (\|\mathbf{x}(k)-\mathbf{x}(k-1)\|,\|\mathbf{x}(k)-\mathbf{x}_{\text{p}}(k)\|)^T,$ (c.f. \eqref{eq:grad_fenchel}),
$\mathbf{z}^{\perp}(k) \triangleq \boldsymbol{\Pi}_{\mathbf{1}_n}^{\perp}\mathbf{z}(k),$ $\Delta_{\mathbf{x}}(k) \triangleq \mathbf{x}(k) - \mathbf{x}_p(k),$ 
and $\boldsymbol{\epsilon}(k) \triangleq \boldsymbol{\Pi}_{\mathbf{1}_n}^{\perp}(\mathbf{x}(k)-\mathbf{x}_p(k)-\mathbf{z}(k)).$
For consistency in the proofs we will define $\mathbf{y}(-1)= \mathbf{y}(0)$ and $\mathbf{x}(-1) = \mathbf{x}(0) = \mathbf{z}(0).$

Note that in Fig. \ref{fig:circle} we have that two arrows stem from $\mathbf{s}$ and end at $\boldsymbol{\epsilon}.$ This is take care by the small gain theorem by accounting for the additive gain of the two paths rather than both of them separately. 

Before starting to prove that each of the arrows holds true we require some intermediate results. Proof of the two lemmas introduced below can be found in \cite{panda}.
\begin{lemma}[Well conditioned Fenchel Conjugate \cite{panda}] The Fenchel conjugate $f^*$ is $\frac{1}{L}-$strongly convex and has $\frac{1}{\mu}$-Lipschitz continuous gradients if and only if the objective function $f$ is $\mu$-strongly convex and has $L-$Lipschitz continuous gradients. 
\end{lemma}
\begin{proof}
See the proof of Theorem 3 in \cite{panda}.
\end{proof}
\begin{lemma}[Equivalent iterates \cite{panda}] The Eco-PANDA iterates can be equivalently written as
\begin{subequations}
\label{eq:iterates}
\begin{align}
\label{eq:iterates:x}
&\mathbf{x}(k+1) := \mathbf{x}(k) - \frac{1}{\eta}(\nabla f(\mathbf{x})-\mathbf{y}(k)) \\
&\mathbf{z}(k+1) := (\mathbf{W}(k)\otimes \mathbf{I}_p)\mathbf{z}(k) + \mathbf{x}(k+1) - \mathbf{x}(k) \\
&\mathbf{y}(k+1) := \mathbf{y}(k)-c(\boldsymbol{\Pi}_{\mathbf{1}_n}^{\perp}\otimes \mathbf{I}_p)(\mathbf{x}(k+1)-\mathbf{z}(k+1))
\end{align}
\end{subequations}
\end{lemma}
\begin{proof}
See the proof of Lemma 2 in \cite{panda}.
\end{proof}
Where the iterates in $\mathbf{x}$ can be seen as a step in the direction of the gradient of \eqref{eq:grad_fenchel} with the time varying $\mathbf{y}(k),$ the iterate in $\mathbf{z}$ tracks the moving average of $\mathbf{x}$ and the iterate in $\mathbf{y}$ corresponds to an inexact projected gradient step.
We now formally express each of the arrows in the form of \eqref{eq:mod}.
\begin{lemma}[$\mathbf{r} \to \mathbf{x}^{\perp}$\label{lemma:a1}]
If $f$ is $\mu-$strongly convex (c.f. Assumption \ref{assumption:function}), it holds that
\begin{equation}
\|\mathbf{x}^{\perp}\|^{\lambda,K} \leq \frac{2\lambda + (1+q)}{\lambda^2 \mu(1-\frac{q}{\lambda})}\|\mathbf{r}\|^{\lambda,K},
\end{equation}
for all $\lambda \in (q,1)$ and $K \geq 0,$ where $q \triangleq 1-\frac{\mu}{\eta}$ for $\eta >L.$
\end{lemma}
The proof of Lemma \ref{lemma:a1} uses the fact that the iterate \eqref{eq:iterates:x} of Eco-PANDA can be seen as running gradient descent on the time varying objective function $f(\mathbf{x}) - \mathbf{y}(k)^T\mathbf{x}.$ In other words, every time the objective function changes, via the change of $\mathbf{y}(k)$ a single gradient step is taken using the previous iteration as a warm start. This strategy has been used to track the time varying optimizer of a sequence of time varying functions in \cite{popkov}. In particular, we use the results in \cite{popkov} to establish a connection between $\mathbf{x}(k+1)$ and $\mathbf{x}_p(k+1)$ (c.f. \eqref{eq:grad_fenchel} for the definition of $\mathbf{x}_p(k)$). Then, we use the results in \cite{panda} to finalize the proof.
\begin{proof}
From Lemma 3 in \cite{panda} we know that
$\|\boldsymbol{\Pi}_{\mathbf{1}_n}^{\perp}\mathbf{x}(k+1)\| \leq \frac{\|\mathbf{y}(k) - \mathbf{y}^{\star}\|}{\mu}.$ This follows from the fact that $\mathbf{x}(k+1) = \nabla f^*(\mathbf{y}(k)),$ $\mathbf{x}^{\star} = \nabla f^*(\mathbf{y}^{\star}) $ and $f^*$ is $\frac{1}{\mu}-$Lipschitz continuous.

 Hence, as an intermediate step to establishing the arrow, we will study the quantity $\|\mathbf{x}(k+1)-\mathbf{x}_p(k+1)\|.$ The iterates in \eqref{eq:iterates:x} can be interpreted as using gradient descent to solve the sequence of optimization problems parametrized by the sequence $\{\mathbf{y}(k)\}_{k \geq 0}.$ More specifically, at each time $k+1$ we perform a step in the direction of the gradient of $f(\mathbf{x})-\mathbf{y}(k)^T\mathbf{x}.$  Therefore, the analysis of how far away $\mathbf{x}(k+1)$ is from $\mathbf{x}_p(k+1)$ can be carried out in a similar fashion as that in \cite{popkov}.
We now adopt the methodology used in the proof of Theorem 1 in \cite{popkov}.

By the fundamental theorem of calculus it follows that
\begin{align}
\label{eq:calculus}
&\nabla f(\mathbf{x}(k))-\mathbf{y}(k) = \nabla f(\mathbf{x}_p(k+1))-\mathbf{y}(k) + \\
&\int_0^1\nabla^2f(\mathbf{x}_p(k+1) + \tau (\mathbf{x}(k)-\mathbf{x}^{\star}(k+1)))\nonumber\\
& \times(\mathbf{x}(k) - \mathbf{x}_p(k+1))\,d\tau \nonumber
\end{align}
 where $\nabla f(\mathbf{x}_p(k+1))-\mathbf{y}(k) = 0$ since $\mathbf{x}_p(k+1)$ is the minimizer of \eqref{eq:grad_fenchel}.
 Let 
 \begin{align}
 \label{eq:amatrix}
 \mathbf{A}_k \triangleq \int_0^1\nabla^2f(\mathbf{x}_p(k+1) + \tau (\mathbf{x}(k)-\mathbf{x}_p(k+1)))\,d\tau.
 \end{align}
 Then, $\mu \mathbf{I}_{np} \preceq \mathbf{A}_k \preceq L\mathbf{I}_{np}$,
 where $\mathbf{I}_{np}$ denotes the size $np$ identity matrix. Therefore, by combining \eqref{eq:calculus} and \eqref{eq:amatrix} we obtain
 \begin{equation}
 \mathbf{A}_k(\mathbf{x}(k) - \mathbf{x}_p(k+1)) = \nabla f(\mathbf{x}(k))-\mathbf{y}(k).
 \end{equation}
Then, we have that
\begin{align}
&\|\boldsymbol{\Pi}_{\mathbf{1}_n}^{\perp}(\mathbf{x}(k+1)-\mathbf{x}_p(k+1))\|= \\
&\nonumber \|(\boldsymbol{\Pi}_{\mathbf{1}_n}^{\perp})(\mathbf{x}(k) - \mathbf{x}_p(k+1)-\frac{1}{\eta}\mathbf{A}_k(\mathbf{x}(k)-\mathbf{x}_p(k+1)))\| \leq \\
& \nonumber\|\mathbf{I}_{np} - \frac{1}{\eta}\mathbf{A}_k\|_{\text{F}}\|\|\boldsymbol{\Pi}_{\mathbf{1}_n}^{\perp}(\mathbf{x}(k)-\mathbf{x}_p(k+1))\|
\end{align}
where $\sigma_{\max}\{\mathbf{I}_{np}-\mathbf{A}_k\} \leq \max \{|1-\frac{\mu}{\eta}|,|1-\frac{L}{\eta}|\}.$ If $\eta > L$ then, the largest of the two elements is $q = 1 - \frac{\mu}{\eta}.$
By using the triangular inequality we can therefore establish that
\begin{align}
\label{eq:tracking}
&\|\boldsymbol{\Pi}_{\mathbf{1}_n}^{\perp}(\mathbf{x}(k+1)-\mathbf{x}_p(k+1))\| \leq \\ &q\left( \|\boldsymbol{\Pi}_{\mathbf{1}_n}^{\perp}(\mathbf{x}(k)-\mathbf{x}_p(k))\|+ \|\boldsymbol{\Pi}_{\mathbf{1}_n}^{\perp}\mathbf{x}_p(k)-\mathbf{x}_p(k+1)\|\right). \nonumber
\end{align}
By further using the triangular inequality we can establish the following statements
\begin{align}
&\|\boldsymbol{\Pi}_{\mathbf{1}_n}^{\perp}\mathbf{x}(k+1)\| - \|\boldsymbol{\Pi}_{\mathbf{1}_n}^{\perp}\mathbf{x}_p(k+1)\| \leq \label{eq:inequality1},\\
& \quad \quad \quad\|\boldsymbol{\Pi}_{\mathbf{1}_n}^{\perp}(\mathbf{x}(k+1)-\mathbf{x}_p(k+1))\| \nonumber \\
& q\|\boldsymbol{\Pi}_{\mathbf{1}_n}^{\perp}(\mathbf{x}(k)-\mathbf{x}_p(k))\| \leq q\|\mathbf{x}^{\perp}(k)\| + q\|\boldsymbol{\Pi}_{\mathbf{1}_n}^{\perp}\mathbf{x}_p(k)\|, \label{eq:inequality2} \\
& \|\boldsymbol{\Pi}_{\mathbf{1}_n}^{\perp}(\mathbf{x}_p(k+1)-\mathbf{x}_p(k))\| \leq  \|\boldsymbol{\Pi}_{\mathbf{1}_n}^{\perp}\mathbf{x}_p(k+1)\| \nonumber\\ & \quad \quad \quad +\|\boldsymbol{\Pi}_{\mathbf{1}_n}^{\perp}\mathbf{x}_p(k)\| \label{eq:inequality3}.
\end{align}
Combining \eqref{eq:tracking}-\eqref{eq:inequality3} we obtain
\begin{align}
&\|\mathbf{x}^{\perp}(k+1)\| \leq q\left(\|\mathbf{x}^{\perp}(k)\|+\|\boldsymbol{\Pi}_{\mathbf{1}_n}^{\perp}\mathbf{x}_p(k)\|\right) + \\
&\|\boldsymbol{\Pi}_{\mathbf{1}_n}^{\perp}\mathbf{x}_p(k+1)\| + \|\boldsymbol{\Pi}_{\mathbf{1}_n}^{\perp}\mathbf{x}_p(k+1))|+\|\boldsymbol{\Pi}_{\mathbf{1}_n}^{\perp}\mathbf{x}_p(k)\|
\end{align}
and therefore
\begin{equation}
\|\mathbf{x}^{\perp}(k+1)\| \leq q\|\mathbf{x}^{\perp}(k)\| + \frac{1+q}{\mu}\|\mathbf{r}^{\perp}(k-1)\| + \frac{2}{\mu}\|\mathbf{r}^{\perp}(k)\|.
\end{equation}
By multiplying both sides by $\lambda^{-(k+1)}$ and taking $\underset{k=1,\hdots,K-1}{\sup}$ we obtain the desired bound for $\lambda \in (q,1).$
\end{proof}
\begin{lemma}[$\mathbf{x}^{\perp} \to \Delta \mathbf{xz}^{\perp}$] Under Assumption \ref{assumption:network} it holds that
\begin{align}
\|\Delta_{xz}^{\perp}\| \leq &\frac{2(1-\lambda^B)}{(1-\lambda)(\lambda^B-\delta)}\|\mathbf{x}^{\perp}\|^{\lambda,K} \nonumber\\
& \frac{\lambda^B}{\lambda^B - \delta}\sum_{t=1}^B \lambda^{1-t}\|\Delta_{xz}^{\perp}(t-1)\|,
\end{align}
for all $\lambda \in (\delta^{\frac{1}{B}},1)$ and for all $K \geq 0.$
\end{lemma}
\begin{proof}
See the proof of Lemma 4 in \cite{panda}.
\end{proof}
\begin{lemma}[$\Delta_{\mathbf{xz}}^{\perp} \to \Delta \mathbf{y}$]
It holds that
\begin{equation}
\|\Delta_{\mathbf{xz}}^{\perp}\|^{\lambda,K} \leq c \|\Delta \mathbf{y}\|^{\lambda,K},
\end{equation}
for any $c >0, \, K>0$ and $\lambda \in (0,1).$
\end{lemma}
\begin{proof}
See Lemma 5 in \cite{panda}.
\end{proof}
\begin{lemma}[\label{lemma:s}$\Delta \mathbf{y} \to \mathbf{s}$]
Given that the objective function $f$ is $\mu-$strongly convex (c.f. Assumption \ref{assumption:function}) it holds that
\begin{equation}
\|\mathbf{s}\|^{\lambda,K} \leq \frac{\lambda \max \{\frac{1}{\eta},\frac{q}{\mu}\}}{\lambda-\rho}\|\Delta \mathbf{y}\|^{\lambda,K}, \text{ for } \lambda \in (\rho,1).
\end{equation}
\end{lemma}
In order to establish how well the average tracking sequence $\mathbf{z}(k+1)$ is performing, we require to understand how the quantity $\mathbf{x}(k+1)-\mathbf{x}(k)$ is behaving. In the case of PANDA, this is fully characterized by the quantity $\mathbf{y}(k)-\mathbf{y}(k-1)$ through Lipschitz continuity of the gradients of the Fenchel conjugate. However, the variation $\mathbf{x}(k+1)-\mathbf{x}(k)$ in the case of Eco-PANDA does not exclusively depend on the difference $\mathbf{y}(k)-\mathbf{y}(k-1)$ but also on how close to $\mathbf{x}_p(k)$ $\mathbf{x}(k)$ was. This is captured in the statement in Lemma \ref{lemma:s}.
\begin{proof}
From the iterates \eqref{eq:iterates:x} it follows that
\begin{equation}
\mathbf{x}(k+1) - \mathbf{x}(k) = -\frac{1}{\eta}\nabla f(\mathbf{x}(k)) + \frac{1}{\eta}\mathbf{y}(k)
\end{equation}
and hence
\begin{equation}
\|\mathbf{x}(k+1)-\mathbf{x}(k)\| \leq \frac{1}{\eta}\|\nabla f(\mathbf{x}(k))-\mathbf{y}(k-1)\|+\frac{1}{\eta}\|\Delta \mathbf{y}(k)\|.
\end{equation}
Note that from \eqref{eq:grad_fenchel} it follows that $\nabla f(\mathbf{x}_p(k))=\mathbf{y}(k-1)$ and hence, by using that $f$ has $L-$Lipschitz continuous gradient, we have 
\begin{equation}
\label{eq:sabove}
\|\mathbf{x}(k+1)-\mathbf{x}(k)\| \leq  \frac{L}{\eta}\|\mathbf{x}(k) - \mathbf{x}_p(k)\| + \frac{1}{\eta} \|\Delta \mathbf{y}(k)\|.
\end{equation}
Further, by using the proof of Theorem 1 in \cite{popkov} and the fact that $f^*$ is $\frac{1}{m}-$Lipschitz continuous, it follows that
\begin{equation}
\label{eq:sabove2}
\|\mathbf{x}(k+1)-\mathbf{x}_p(k+1)\| \leq q\left(\|\mathbf{x}(k) - \mathbf{x}_p(k)\| + \frac{1}{m}\|\Delta \mathbf{y}(k)\|\right).
\end{equation}
The inequalities in \eqref{eq:sabove} and \eqref{eq:sabove2} can be expressed jointly as 
\begin{equation}
\label{eq:matrix}
\mathbf{s}(k+1) \leq \begin{pmatrix} 0 & \frac{L}{\eta} \\
0 & q 
\end{pmatrix}
\mathbf{s}(k) + \max \left\{\frac{1}{\eta},\frac{q}{m}\right\}\|\Delta \mathbf{y}(k)\|.
\end{equation}
Note that the spectral radius of the matrix appearing in \eqref{eq:matrix} is upper bounded by $\rho \triangleq \max\{\frac{L}{\eta},q\}$ and hence it holds that
\begin{equation}
\|\mathbf{s}(k+1)\| \leq \rho \|\mathbf{s}(k)\| + \max \left\{\frac{1}{\eta},\frac{q}{\mu}\right\}.
\end{equation}
By performing the same steps as described in the proof of Lemma \ref{lemma:a1} we obtain the desired result.
\end{proof}
\begin{lemma}[$\mathbf{s} \to \Delta_{\mathbf{x}}$]
It holds that $\|\Delta_{\mathbf{x}}\|^{\lambda,K} \leq \|\mathbf{s}\|^{\lambda,K}.$
\end{lemma}
\begin{proof}
The bound follows from the definition of $\mathbf{s}.$
\end{proof}
\begin{lemma}[$\mathbf{s} \to \mathbf{z}^{\perp}$]
Under assumption \ref{assumption:network} it holds that
\begin{align}
\|\mathbf{z}^{\perp}\|^{\lambda,K} \leq \frac{1-\lambda^B}{(1-\lambda)(\lambda^B-\delta)}\|\mathbf{s}\|^{\lambda,K}  \nonumber \\
+ \frac{\lambda^B}{\lambda^B - \delta}\sum_{t=1}^B \lambda^{1-t}\|\mathbf{z}(t-1)\|,
\end{align}
for $\lambda^B \in (\delta,1).$
\end{lemma}
\begin{proof}
See the proof of Lemma 6 in \cite{panda}.
\end{proof}
\begin{lemma}[$\Delta_{\mathbf{x}},\,\mathbf{z}^{\perp} \to \boldsymbol{\epsilon}$]
It holds that $\|\boldsymbol{\epsilon}\|^{\lambda,K} \leq \|\Delta_{\mathbf{x}}\|^{\lambda,K} + \|\mathbf{z}^{\perp}\|^{\lambda,K}.$
\end{lemma}
\begin{proof}
The proof follows by using the triangle inequality on the definition of $\boldsymbol{\epsilon}.$
\end{proof}
\begin{lemma}[$\boldsymbol{\epsilon} \to \mathbf{r}$]
Under Assumptions \ref{assumption:network} and \ref{assumption:function}, $c \in \left( 0,\frac{\mu}{2}\right]$ and $\lambda \in \left[\sqrt{1-\frac{c}{2L}},1 \right)$ it holds that
\begin{equation}
\|\mathbf{r}\|^{\lambda,K} \leq \sqrt{L\mu}\|\boldsymbol{\epsilon}\|^{\lambda,K} + 2 \|\mathbf{r}(0)\|.
\end{equation}
\end{lemma}
\begin{proof}
See Lemma 8 in \cite{panda}.
\end{proof}
We now have established all the arrows in Fig. \ref{fig:circle} and are ready to invoke the small gain theorem. In order to be able to claim that Eco-PANDA converges R-linearly we have to establish that there exists a tuple $(\lambda,c,\eta)$ such that the conditions
\begin{align}
\label{eq:enormous}
&\frac{2(2\lambda+1+q)(1-\lambda^B)A(1-\lambda^B + (1-\lambda)(\lambda^B-\delta))}{(\lambda-q)(\lambda-\rho)((1-\lambda)(\lambda^B-\delta))^2} \nonumber\\
&\quad\quad \times c < \frac{1}{\sqrt{\kappa}}, \text{ with } A \triangleq \max\left\{\frac{1}{\eta},\frac{q}{\mu}\right\}\\
& \lambda \in \left(\max\left\{q,\rho,\delta^{1/B},\sqrt{1-\frac{c}{2L}}\right\},1\right) \label{eq:lambda_interval}, \quad c \in \left(0,\frac{\mu}{2}\right]
\end{align}
for $q = 1-\frac{\mu}{\eta},\, \eta >0 $ are fulfilled.

By using that for $0.5 \leq \lambda \leq 1$ it holds that $(1-\lambda^B)/(1-\lambda) \leq B$ \cite{DIGing}, we have that if $c$ fulfils 
\begin{equation}
c \leq \frac{(\lambda-q)(\lambda-\rho)(\lambda^B-\delta)^2}{\sqrt{\kappa}(A + B(\lambda^B - \delta))(2(2\lambda+1+q))}
\end{equation}
it will also fulfil \eqref{eq:enormous}. Recall now that $q = 1 -\frac{\mu}{\eta}$ and $\rho = \max\{q,\frac{L}{\eta}\}.$ Further, by inspecting \eqref{eq:lambda_interval} we may conclude that $\lambda \geq \sqrt{1-\frac{\mu}{4L}} \triangleq \lambda_{\min}.$ In particular, by choosing $\eta = 4L$ we have that $\sqrt{q} = \lambda_{\min}$ and also $\rho = q.$ Therefore, we conclude that if we chose $\eta = 4L$ we have that any $c$ fulfilling
\begin{equation}
c \leq \frac{(\sqrt{q}-q)^2(\lambda^B-\delta)^2}{2\sqrt{\kappa}(A+B(\lambda^B-\delta))(3+q)}
\end{equation}
fulfils \eqref{eq:enormous}.
Note that we use the fact that the term $(\lambda-q)^2$ is monotonically increasing for $\lambda \geq \sqrt{q}.$ Further, note that as long as the quantity $\eta > L$ we have that $A < \frac{1}{\mu}$ and therefore
\begin{equation}
c \leq \frac{\mu(\sqrt{q}-q)^2(\lambda^B-\delta)^2}{2\sqrt{\kappa}(\mu A + B\mu (1 - \delta))(3+q)}.
\end{equation}
It is straight forward to show that the quantity $\frac{(\sqrt{q}-q)^2}{q(3+q)}<1$ for all non-negative values of $q$ and by recalling that $A \geq \frac{q}{\mu}$ we may conclude that if $c$ fulfils
\begin{equation}
c \leq \frac{\alpha \mu(\lambda^B-\delta)^2}{2\sqrt{\kappa}}
\end{equation}
for an appropriately chosen $\alpha \in (0,1),$ it also fulfils \eqref{eq:enormous}.
This implies that we are to establish that there exists a pair $(\lambda,c)$ fulfilling
\begin{align}
c \leq \frac{\alpha \mu(\lambda^B - \delta)^2}{2\sqrt{\kappa}},\quad
c \geq 2L(1-\lambda^2) \\
\lambda \in \left(\delta^{1/B},1 \right), \quad
c \in \left(0,\frac{\mu}{2} \right] 
\end{align}
given that $\eta = 4L.$ Then, the procedure to establish the result is identical to that in Section 4.C in \cite{panda} but taking into account the parameter $\alpha.$
\section{Numerical Experiments}
In order to illustrate the performance of Eco-PANDA, consider the decentralized estimation problem
\begin{equation}
\label{eq:numerical}
\underset{\mathbf{x} \in \mathbb{R}^p}{\min} \quad \frac{1}{2nd}\sum_{i=1}^n \|\mathbf{H}_i\mathbf{x}-\mathbf{b}_i\|^2+ \frac{r}{2}\|\mathbf{x}\|^2.
\end{equation}
The elements of the vector $\mathbf{\hat{x}}\in \mathbb{R}^d$ are generated independently and according to the distribution $\mathcal{N}(0,10).$ The elements of the matrices $\mathbf{H}_i$ are generated independently and according to the distribution $\mathcal{N}(0,0.1).$ Then, the vectors $\mathbf{b}_i$ are generated according to the equation $\mathbf{b}_i = \mathbf{H}_i\mathbf{\hat{x}} + \boldsymbol{\xi}_i$ where the elements of $\boldsymbol{\xi}_i$ are generated independently and according to the distribution $\mathcal{N}(0,0.1).$ The sequence of graphs $\mathcal{G}(\mathcal{V},\mathcal{E}(k))$ on which the problem \eqref{eq:numerical} is solved is randomly generated. In particular, each bi-directional arc will belong to the set of edges with probability $\pi.$ The Metropolis Hastings matrix is chosen for both methods, i.e. the mixing matrix is given by
\begin{equation}
w_{i,j}(k) = \begin{cases}
\frac{1}{\max\{d_i(k),d_j(k)\}}& (i,j) \in \mathcal{E}(k) \\
1-\sum_{j\in \mathcal{N}_i(k)}w_{i,j}(k) & i = j \\
0 & \text{otherwise.}
\end{cases}
\end{equation}
\begin{figure}
\begin{center}
\scalebox{0.53}{\input{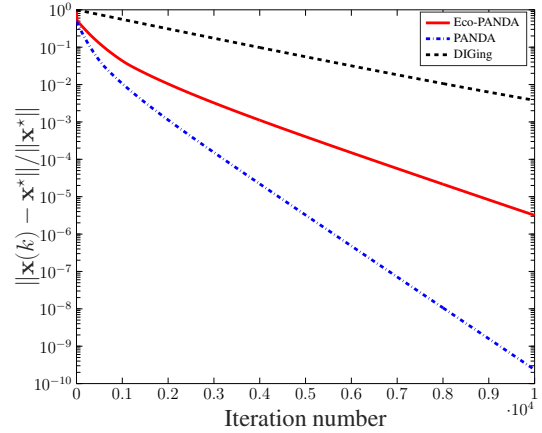}}
\end{center}
\caption{\label{fig:plot} Performance of Eco-PANDA, PANDA and DIGing.}
\end{figure}
Fig. \ref{fig:plot} is generated by setting $n = 10,$ $p = 3,$ $d = 5,$ $r = 1$ and $\pi = 0.1.$ Picking $r = 1$ allows for the optimization problems to be well conditioned allowing for PANDA and Eco-PANDA to perform well. On ill conditioned problems, DIGing is expected to perform better \cite{panda}. This observation is consistent with the empirical evidence in \cite{panda} and the analysis in both \cite{panda} and \cite{DIGing}, where in \cite{panda} we require the strong convexity of each of the function components and in \cite{DIGing} the strong convexity of only one of them is required.
 The step-sizes for the methods have been hand-optimized and are set to $c = 0.0005,$ $\eta = 0.5$ for Eco-PANDA, $c_{\text{PANDA}} = 0.001$ for PANDA and $\alpha_{\text{DIGIing}} = 0.003$ for DIGing respectively.

From Fig. \ref{fig:plot} one can see that Eco-PANDA is outperformed by PANDA. This is reasonable as the iterates of Eco-PANDA are an approximation of those of PANDA. Note however that Eco-PANDA could still be preferable in terms of computation cost due to the simpler iterates. Further, as could be concluded from Theorem \ref{theorem:main} Eco-PANDA tolerates a smaller step size before starting to diverge than PANDA. However, the iterates of PANDA require a matrix inversion at each iterate to solve the optimization problem \eqref{eq:fencheldual} while the iterates of Eco-PANDA require just the computation of the gradient.
\label{sec:typestyle}
\section{Conclusions}
In this paper we proposed Eco-PANDA a linearized version of PANDA that we have proven converges linearly on time-varying graphs. By doing this we have circumvented the main disadvantage PANDA has when compared to DIGing, that is, that the iterates are computationally speaking more expensive. Eco-PANDA alleviates the nodes' computational burden while keeping the low communication cost of PANDA.  We analytically and numerically observe the performance loss due to the having simpler iterates. The loss in performance is directly related to the shrinking of the largest allowed step-size which provides for the best convergence bound.
\bibliographystyle{IEEEbib}

\end{document}